\documentclass[a4paper]{amsart}

\title[Representations Theory of RCA via Microlocal Analysis]{%
Representation theory of the rational Cherednik algebras of type $\Z/l\Z$
via microlocal analysis}
\author{Toshiro Kuwabara}
\thanks{The author is partially supported by 
Grant-in-Aid for Young Scientist (B) 21740013,
Japan Society for the Promotion of Science.}
\address{Research Institute for Mathematical Science, 
Kyoto University, \\
Oiwakecho Kitashirakawa, Sakyoku, Kyoto, 606-8502, JAPAN}
\email{toshiro@kurims.kyoto-u.ac.jp}

\usepackage{latexsym}
\usepackage{amsmath}
\usepackage{amsfonts}
\usepackage{amssymb}
\usepackage{amsxtra}
\usepackage{mathrsfs}
\usepackage [matrix,arrow]{xy}

\newtheorem{definition}{Definition}[section]
\newtheorem{proposition}[definition]{Proposition}
\newtheorem{theorem}[definition]{Theorem}
\newtheorem{corollary}[definition]{Corollary}
\newtheorem{lemma}[definition]{Lemma}
\newtheorem{remark}[definition]{Remark}

\newcommand{\refdef}[1]{Definition~\ref{#1}}
\newcommand{\refprop}[1]{Proposition~\ref{#1}}
\newcommand{\refthm}[1]{Theorem~\ref{#1}}
\newcommand{\refcor}[1]{Corollary~\ref{#1}}
\newcommand{\reflemma}[1]{Lemma~\ref{#1}}
\newcommand{\refremark}[1]{Remark~\ref{#1}}
\newcommand{\refeq}[1]{(\ref{#1})}
\newcommand{\refsec}[1]{Section~\ref{#1}}

\newcommand{\refappendix}[1]{Appendix~\ref{#1}}

\newcommand{\defeq}{\ensuremath{\underset{\mathrm{def}}{=}}}
\newcommand{\C}{{\mathbb C}}

\newcommand{\Z}{{\mathbb Z}}

\newcommand{\Dim}{\mathrm{dim}}

\DeclareMathOperator{\Spec}{Spec}
\DeclareMathOperator{\Res}{Res}

\newcommand{\bbT}{{\mathbb T}}

\newcommand{\frS}{\mathfrak{S}}

\newcommand{\g}{\mathfrak{g}}
\newcommand{\frh}{\mathfrak{h}}

\newcommand{\ii}{\eta}

\newcommand{\bbP}{{\mathbb P}}
\newcommand{\scL}{\mathscr{L}}
\newcommand{\scM}{\mathscr{M}}
\newcommand{\scN}{\mathscr{N}}
\newcommand{\A}{\mathscr{A}}
\newcommand{\tildA}{\widetilde{\mathscr{A}}}
\newcommand{\W}{\mathscr{W}}
\newcommand{\calD}{\mathcal{D}}
\newcommand{\calO}{\mathcal{O}}
\newcommand{\calL}{\mathcal{L}}
\newcommand{\calM}{\mathcal{M}}
\newcommand{\calN}{\mathcal{N}}

\newcommand{\bfe}{\mathbf{e}}
\newcommand{\bfDelta}{\mathbf{\Delta}}

\newcommand{\bbG}{{\mathbb G}}
\newcommand{\scF}{\mathscr{F}}
\DeclareMathOperator{\Supp}{Supp}

\DeclareMathOperator{\gEnd}{End}
\DeclareMathOperator{\lEnd}{\operatorname{\mathscr{E}\kern-.1pc\mathit{nd}}}
\DeclareMathOperator{\gHom}{Hom}
\DeclareMathOperator{\lHom}{\operatorname{\mathscr{H}\kern-.1pc\mathit{om}}}
\DeclareMathOperator{\Mod}{Mod}
\newcommand{\mmod}{\text{-}\mathrm{mod}}
\DeclareMathOperator{\Ad}{Ad}
\DeclareMathOperator{\ad}{ad}
\DeclareMathOperator{\Aut}{Aut}
\DeclareMathOperator{\Der}{Der}
\DeclareMathOperator{\Irr}{Irr}

\newcommand{\opp}{\mathrm{opp}}

\newcommand{\isoto}[1][]{\mathop{\xrightarrow[#1]%
{\rule{0pt}{.9ex}%
{\raisebox{-.4ex}[0ex][-.6ex]{$\mspace{3mu}\sim\mspace{3mu}$}}}}}
\newcommand{\scbul}{{\,\raise1pt\hbox{$\scriptscriptstyle\bullet$}\,}}
\newcommand{\prolim}[1][]{\mathop{\varprojlim}\limits_{#1}}

\newcommand{\One}{\mathbf{1}}

\begin{document}
\begin{abstract}
 Based on the methods developed in \cite{KR},
 we consider microlocalization of the rational Cherednik algebra of
 type $\Z/l\Z$. Our goal is to construct 
 the irreducible modules and standard modules of the rational Cherednik
 algebra by
 using the microlocalization. As a consequence, we obtain the 
 sheaves corresponding to holonomic systems with
 regular singularities.
\end{abstract}

\maketitle

\section{Introduction}

The symplectic reflection algebra is a noncommutative deformation of
the smash product $\C[V] \# \Gamma$, introduced by \cite{EG}, where
$V$ is a symplectic vector space and $\Gamma$ is a finite group 
generated by symplectic reflections of $V$. Sometimes, we identify
the symplectic reflection algebra with its spherical subalgebra, a
noncommutative deformation of $\C[V]^{\Gamma}$, because these algebras
are mutually Morita equivalent except for a certain choice of their 
parameters.

When the group $\Gamma$ coincides with a complex reflection group 
and $V$ coincides with $\frh \oplus \frh^{*}$ where $\frh$ is the
reflection representation of $\Gamma$, the symplectic reflection algebra
is sometimes called the rational Cherednik algebra. An important property
of the rational Cherednik algebra is that it has a triangular decomposition 
similarly to complex semisimple Lie algebras. Via the triangular decomposition,
we can introduce a certain subcategory of the category of modules,
called the category $\calO$. The category $\calO$ is a highest weight
category in terms of \cite{CPS}. 
Its standard modules and costandard 
modules are studied in \cite{GGOR}. 
For each irreducible $\C \Gamma$-module $E \in \Irr \C\Gamma$, 
we have a corresponding 
standard module $\Delta(E)$. The standard module has a unique irreducible 
quotient $L(E)$ and any irreducible module in the category $\calO$ is 
isomorphic to $L(E)$ for a certain $E$. 
One of fundamental problems of the representation theory of the 
rational Cherednik algebras is to determine multiplicities $[\Delta(E) : L(F)]$ 
 in the Grothendieck group for $E$, $F \in \Irr \C\Gamma$.

When the group $\Gamma$ is a wreath product $(\Z/l\Z) \wr \frS_n$ of a
cyclic group $\Z/l\Z$ and a symmetric group $\frS_n$, there is a close
connection between the rational Cherednik algebra and a quiver variety
which is a symplectic variety introduced by \cite{Na}. 
After the leading work of \cite{GS1} and \cite{GS2}, \cite{KR} constructed
a microlocalization of the rational Cherednik algebra of type $\frS_n$.
The microlocalization is a kind of Deformation-Quantization algebra,
called a W-algebra, on a quiver variety. \cite{KR} introduced the
notion of F-action of the W-algebras and established the equivalence
of categories between the category of modules of the rational Cherednik
algebra and the category of F-equivariant, good modules of the W-algebra.
This equivalence is an analogue of the Beilinson-Bernstein correspondence
for complex semisimple Lie algebras.

In \cite{BK}, microlocalization of the rational Cherednik algebra of 
type $\Z/l\Z$ was studied. As an application of the microlocalization 
of the rational Cherednik algebras, we study the construction of
the irreducible modules and the standard modules of the rational Cherednik algebra
via microlocalization.

Let us describe the structure of this article.

In \refsec{sec:quiver-varieties}, we review fundamental properties
of the minimal resolutions
of Kleinian singularities of type A. We construct the Kleinian
singularities and their resolutions $X$ as quiver varieties of a cyclic
quiver. Moreover, we see that the structure of $X$
as a toric variety gives us an affine open covering 
$X = \bigcup_{i=1}^{l} X_i$ such that $X_i \simeq \C^2$.

In \refsec{sec:w-algebra}, we review the general setting of 
the W-algebra and construct the microlocalization $\tildA_c$ of the rational
Cherednik algebra $A_c$ on $X$.
By the theorem of \cite{BK}, we have an equivalence of categories
\begin{align*}
 \Mod_{F}^{good}(\tildA_c) &\longrightarrow A_c\mmod \\
 \scM &\mapsto \gHom_{\Mod_F^{good}(\tildA_c)}(\tildA_c, \scM)
\end{align*}
under certain conditions on the parameter $c$.

At the end of \refsec{sec:w-algebra-on-quiver-var}, we describe the
structure of $\tildA_c \vert_{X_i}$ explicitly 
on the affine open subset $X_i$ for $i=1$, $\dots$, $l$.

In \refsec{sec:RCA}, we briefly review the representation 
theory of the rational Cherednik algebra. Its spherical subalgebra
is isomorphic to $A_c$. 
We introduce the category $\calO(A_c)$, and review
the definition of its standard modules $\Delta_c(i)$ and 
irreducible modules $L_c(i)$. 

In \refsec{sec:constr-stand-mod}, we construct an $\tildA_c$-module
$\calM^{\Delta}_c(i)$ for $i=1$, $\dots$, $l$. It is also an F-equivariant,
holonomic $\tildA_c$-module supported on a certain Lagrangian 
subvariety. 
We show that the corresponding $A_c$-module 
$\gHom_{\Mod_F^{good}(\tildA_c)}(\tildA_c, \calM_c^{\Delta}(i))$ is
isomorphic to the standard module $\Delta_c(i)$. 

In \refsec{sec:constr-irred-modul}, we construct an $\tildA_c$-module
$\calL_c(i)$ for $i=1$, $\dots$, $l$. It is an F-equivariant, holonomic
$\tildA_c$-module supported on a certain Lagrangian subvariety. 
Moreover, we show that 
$\calL_c(i)$ is a irreducible $\tildA_c$-module for any $i=1$, $\dots$, $l$.
At the end of \refsec{sec:constr-irred-modul}, we determine the
multiplicity $[\Delta_c(i) : L_c(j)]$ in the Grothendieck group
of $\calO(A_c)$ as a corollary of the construction of
$\tildA_c$-module $\calM_c^{\Delta}(i)$ and $\calL_c(j)$. 

Finally, in \refappendix{sec:glob-sect-modul}, we determine the global
sections of $\tildA_c$-module $\calM_c^{\Delta}(i)$ explicitly. 

\subsection*{Acknowledgments}
The author is deeply grateful to Masaki Kashiwara for his valuable comments and
kind advice about how to manage irreducibility of holonomic
$\W$-modules.
 He thanks Raphael Rouquier for valuable discussions.
He also thanks Tomoyuki Arakawa, Gwyn Bellamy and  Testuji Miwa for valuable discussions
and comments.
The auther was supported by JSPS Grant-in-Aid for Young Scientists (B) 21740013.

\section{Quiver varieties}
\label{sec:quiver-varieties}

In this section we review the definition and fundamental 
properties of quiver varieties without framing, which were 
introduced by \cite{Kr}. 

Let $Q=(I, E)$ be a cyclic quiver with vertices 
$I = \{I_i \;\vert\; i=0, \dots, l-1\}$ and arrows
$E = \{\alpha_i : I_{i-1} \rightarrow I_i \;\vert\; i=1,\dots,l \}$.
Let $\overline{Q} = (I, E \sqcup E^*)$ be a quiver with vertices $I$
and arrows $E$ and $E^* = \{ \alpha^*_i : I_{i} \rightarrow I_{i-1}\}$.
Throughout this paper, we regard indices of vertices and 
edges of $Q$ and $\overline{Q}$ as integers modulo $l$, i.e. we regard
$I_{l+i} = I_i$ and $\alpha_{l+i} = \alpha_i$. 
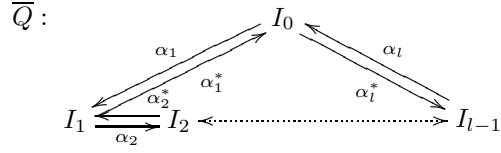
\[
\overline{Q}:
 \xymatrix{
& & I_0 \ar@<-2pt>[dll]_{\alpha_1} \ar@<-2pt>[drr]_{\alpha_l^*} & & \\
I_1 \ar@<-2pt>[r]_{\alpha_2} \ar@<-2pt>[urr]_(0.6){\alpha^*_1}
& I_2 \ar@{<.>}[rrr] \ar@<-2pt>[l]_<{\alpha_2^*} & & &
I_{l-1} \ar@<-2pt>[ull]_{\alpha_l}
}
\]

A representation of $\overline{Q}$ with a dimension vector 
$\delta = (1, \dots, 1)$ is a pair 
($V$, $(a_i, b_i)_{i=1, \dots, l}$)
of an $I$-graded vector space $V = \bigoplus_{i \in 0}^{l-1} V_i$ such
that $\Dim V_i = 1$ for all $i \in I$ and linear maps 
$a_i : V_{i-1} \rightarrow V_i$, $b_i : V_i \rightarrow V_{i-1}$.
Since $\Dim V_i = 1$ for all $i$, we regard $a_i$ and $b_i$ as 
elements of $\C$. 
Let $GL(\delta) = \prod_{i=0}^{l-1} GL(V_i) \simeq (\C^*)^{l}$ 
be a reductive algebraic group
acting on $V$. Let $G = PGL(\delta) = GL(\delta) / \C_{diag}^* \simeq 
(\C^*)^{l-1}$ where
$\C_{diag}^*$ is the diagonal subgroup of $GL(\delta)$. 
Let $\g = \mathrm{Lie}(G)$ be the Lie algebra of $G$.
We have $\g = \left(\bigoplus_{i=0}^{l-1} \C \right) / \C_{diag}$ 
where $\C_{diag}$ is the diagonal of $\bigoplus_{i=0}^{l-1} \C$.

\[
\text{A representation of $\overline{Q}$:}
 
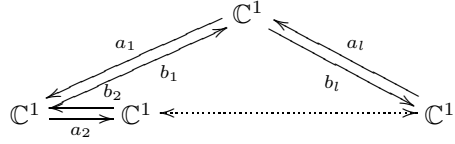
\xymatrix{
& & \C^1 \ar@<-2pt>[dll]_{a_1} \ar@<-2pt>[drr]_{b_l} & & \\
\C^1 \ar@<-2pt>[r]_{a_2} \ar@<-2pt>[urr]_(0.6){b_1}
& \C^1 \ar@{<.>}[rrr] \ar@<-2pt>[l]_<{b_2} & & &
\C^1 \ar@<-2pt>[ull]_{a_l}
}
\]

Fix a parameter $\theta = (\theta_0, \dots, \theta_{l-1}) \in \Z^l$ such
that $\theta_0 + \theta_1 + \dots + \theta_{l-1} = 0$. Note that 
we regard indices of $\theta$ as integers modulo $l$: 
$\theta = (\theta_i)_{i \in \Z/l\Z}$, and 
$\theta_i + \theta_{i+1} + \dots + \theta_{j-1}$ is well-defined
for any $i$, $j \in \Z/l\Z$. We regard 
$\theta$ as an infinitesimal character of $\g$.

A representation $(V, (a_i, b_i)_{i=1, \dots, l})$
is called $\theta$-semistable if any $I$-graded subspace $W$ of $V$
which is stable under the action of $(a_i, b_i)_{i=1, \dots, l}$
satisfies the condition: $\sum_{i=0}^{l-1} \Dim W_i \theta_i \le 0$.

Fix a parameter $\theta$.
Let $\widetilde{X}_\theta \subset \C^{2l}$ 
be the space of all $\theta$-semistable representation,
\[
 \widetilde{X}_\theta = \{(a_i, b_i)_{i=1,\dots,l} \in \C^{2l}
\;\vert\; (a_i, b_i)_{i=1, \dots, l} \text{ is $\theta$-semistable}\}.
\]
The group $G$ acts effectively on $\widetilde{X}$. This action is
a symplectic action. Two points $p$, $p'$ of $X$ is called S-equivalent
if the closures of their orbits have an intersection in 
$\widetilde{X}_\theta$.

Consider the following 
moment map with respect to the action:
\begin{align*}
 \mu: \widetilde{X}_\theta &\longrightarrow \g^* \subset \C^l , \\
 (a_i, b_i)_{i=1, \dots, l} &\mapsto (a_{i+1} b_{i+1} - a_i b_i)_{i=0, \dots, l-1}.
\end{align*}
We consider the Hamiltonian reduction with respect to the moment map
$\mu$. The subset $\mu^{-1}(0) \subset \widetilde{X}_\theta$ is stable under
the action of $G$. 

\begin{definition}
 \label{def:1}
 The quiver variety of the quiver $\overline{Q}$ 
 with the dimension vector $\delta$
 and the stability parameter $\theta$ is a complex symplectic variety
 \[
  X_\theta = \mu^{-1}(0) / \sim_S
 \]
 where $\sim_S$ be the S-equivalence.
\end{definition}

We denote an S-equivalence class in $X_{\theta}$ containing 
$(a_i, b_i)_{i=1, \dots, l} \in \widetilde{X}_\theta$ by
$[a_i, b_i]_{i=1, \dots, l}$.

Let us consider the case of $\theta = 0 = (0, \dots, 0)$. 
For $(a_i, b_i)_{i=1,\dots,l} \in \mu^{-1}(0) \subset \widetilde{X}_0 = \C^{2l}$,
we set $\bar{a} = \sqrt[l]{a_1 \cdots a_l}$, 
$\bar{b} = \sqrt[l]{b_1 \cdots b_l}$ such that 
$\bar{a} \bar{b} = a_1 b_1$. 
Then, we have the following 
isomorphism of algebraic varieties.
\begin{align*}
 X_0 &\xrightarrow{\simeq} \C^2 / (\Z/l\Z) \\
 [a_i, b_i]_{i=1, \dots, l} &\mapsto (\bar{a}, \bar{b})
\end{align*}
Note that the image of $(\bar{a}, \bar{b})$ in $\C^2 / (\Z/l\Z)$ is
independent of the choice of root. 
(cf. \cite{CS}).

\begin{proposition}[\cite{CS}, \cite{Kr}]
\label{prop:1}
 If a stability parameter $\theta = (\theta_i)_{i=0, \dots, l-1}$ satisfies
$\theta_i + \theta_{i+1} + \dots + \theta_{j-1} \ne 0$ for
all $i$, $j$ ($i \ne j$),
$X_\theta$ is
nonsingular and we have a minimal resolution of Kleinian
singularities of type $A_{l-1}$:
\[
\pi_\theta: X_\theta \longrightarrow X_0 \simeq \C^2 / (\Z/l\Z).
\]
\end{proposition}

In the rest of this paper, we fix a stability
parameter $\theta$ satisfying the condition of 
\refprop{prop:1}. We denote 
$X_\theta$ by $X$, $\pi_\theta$ by $\pi$, ...
for simplicity. 

One of the fundamental properties of $X$ is that it is
a toric variety with respect to the following action of a $2$-dimensional
torus $\bbT^2 = (\C^*)^2$:
\[
 (q_1, q_2) [a_i, b_i]_{i=1, \dots, l} = 
 [q_1 a_i, q_2 b_i]_{i=1, \dots, l}
\]
for $(q_1, q_2) \in \bbT^2$ and $[a_i, b_i]_{i=1, \dots, l} \in X$.
The following facts are easy to obtain from the general theory of
toric varieties. Refer to \cite[Section 2]{K} for  proofs of these facts,
or to \cite{Fu} for the general theory of toric varieties. 

The variety $X$ has $l$ $\bbT$-fixed points 
$p'_1$, \dots, $p'_{l}$ where $p'_i = [a_j, b_j]_{j=1, \dots, l}$ 
is given as follows:
\begin{gather}
  a_i = 0,\quad b_i = 0, \nonumber\\
  a_j = 0,\quad b_j \ne 0 \quad \text{if $\theta_{i} + \theta_{i+1} +
 \dots + \theta_{j-1} < 0$}, \nonumber\\
  a_j \ne 0,\quad b_j = 0 \quad \text{if $\theta_{i} + \theta_{i+1} +
 \dots + \theta_{j-1} > 0$}. \nonumber
\end{gather}
Note that we have $\theta_i + \theta_{i+1} +
\dots + \theta_{j-1} \ne 0$ for all $i \ne j$ 
by the condition of the stability parameter $\theta$.

Define an ordering $\rhd$ on the set of indices 
$\Lambda = \{1, \dots, l\}$ by 
\begin{equation*}
 i \rhd j \Longleftrightarrow 
 \theta_{i} + \dots + \theta_{j-1} < 0.
\end{equation*}
By the condition of the stability parameter 
$\theta$, the ordering $\rhd$ is a total ordering. 
Let $\ii_1$, $\dots$, $\ii_l$ be the indices in $\Lambda$ such that
\begin{equation}
\label{eq:2}
 \ii_1 \rhd \ii_2 \rhd \dots \rhd \ii_l.
\end{equation}

\begin{remark}
 Note that the order of the numbering $\ii_1$, $\dots$, $\ii_i$
 is reversed from the one of \cite{K}. 
\end{remark}

Set $p_i = p'_{\ii_i}$ for $i=1$, $\dots$, $l$. The explicit 
description of the point $p_i$ is given as follows.

\begin{lemma}
\label{lemma:1}
 For $i=1$, $\dots$, $l$, the fixed point 
$p_{i} = p'_{\ii_i} = [a_j, b_j]_{j=1, \dots, l}$ is given by
 \begin{gather*}
  a_{\ii_i} = 0, \quad b_{\ii_i} = 0, \\
  a_{\ii_j} = 0, \quad b_{\ii_j} \ne 0 \qquad \text{for $j > i$}, \\
  a_{\ii_j} \ne 0, \quad b_{\ii_j} = 0 \qquad \text{for $j \le i$}.
 \end{gather*}
\end{lemma}

Let us consider a Lagrangian subvariety 
$\pi^{-1}(\{\bar{a}=0 \text{ or } \bar{b}=0\})$. 
This subvariety has $(l+1)$ irreducible components $D_0$, $D_1$, $\dots$,
$D_l$ such that $D_0$, $D_l \simeq \C^1$, $D_i \simeq \bbP^1$
for $1 \le i \le l-1$ and $p_i$ is a unique intersection of
$D_{i-1}$ and $D_i$. We can describe $D_i$ explicitly as follows.

\begin{lemma}
 \label{lemma:2}
 For $i=1$, $\dots$, $l$, the $\bbT$-divisor
 $D_i$ is given by
\[
  D_i = 
 \overline{\left\{ [a_j,b_j]_{j=1,\dots,l} 
 \biggm|
 \begin{array}{l}
  a_{\ii_j} = 0, \; b_{\ii_j} \ne 0 \qquad \text{for $j > i$}, \\
  a_{\ii_j} \ne 0, \; b_{\ii_j} = 0 \qquad \text{for $j \le i$}.
 \end{array}
 \right\}}.
\]
Similarly, $D_0$ is given by
\[
  D_0 = 
 \overline{\left\{ [a_j,b_j]_{j=1,\dots,l} \;\vert\;
 b_j \ne 0
 \right\}}.
\]
\end{lemma}

The description as a toric variety gives us the following affine open 
covering of $X$,
\[
 X = \bigcup_{i=1}^{l} X_i, \qquad
 X_i = 
 \left\{
 [a_j, b_j]_{j=1, \dots, l} \Bigm\vert
 \begin{array}{cc}
  a_{\ii_j} \ne 0 &\text{for } j < i \\
  b_{\ii_j} \ne 0 &\text{for } j > i
 \end{array}
 \right\}.
\]
We introduce coordinate functions
$\bar{x}_i$ (resp. $\bar{y}_i$) for $i=1$, $\dots$, $l$ on $\widetilde{X} \subset \C^{2l}$
defined by $\bar{x}_i((a_j, b_j)_{j=1, \dots, l}) = a_i$ 
(resp. $\bar{y}_i((a_j, b_j)_{j=1, \dots, l}) = b_i$).
For $i=1$, $\dots$, $l$, let $R_i$ be the
following subring of $\C(\bar{x}_1, \dots, \bar{x}_l, \bar{y}_1, \dots, \bar{y}_l)$
which is isomorphic to a polynomial ring in $2$-variables:
\[
 R_i = \C[\bar{f}_i, \bar{g}_i]
\]
where 
\[
 \bar{f}_i = \frac{\bar{x}_{\ii_1} \bar{x}_{\ii_2} \dots \bar{x}_{\ii_i}}
 {\bar{y}_{\ii_{i+1}} \bar{y}_{\ii_{i+2}} \dots \bar{y}_{\ii_l}}, 
 \quad
 \bar{g}_i = \frac{\bar{y}_{\ii_i} \bar{y}_{\ii_{i+1}} \dots \bar{y}_{\ii_l}}
 {\bar{x}_{\ii_{1}} \bar{x}_{\ii_{2}} \dots \bar{x}_{\ii_{i-1}}}
 \in \C(\bar{x}_1, \dots, \bar{x}_l, \bar{y}_1, \dots, \bar{y}_l)
\]
Then we have
\[
 X_i = \Spec R_i \simeq \C^2 = T^* \C^1.
\]
Note that $\bar{f}_i \bar{g}_{i+1} = 1$ on $X_i \cap X_{i+1}$.

For $i=1$, $\dots$, $l$, the fixed point $p_i$ belongs to $X_i$.
For $i=1$, $\dots$, $l-1$, we have $D_i \simeq \bbP^1 \subset X_i \cup X_{i+1}$
and there is an isomorphism 
$X_i \cup X_{i+1} \simeq T^* \bbP^1$. 

\section{W-algebra}
\label{sec:w-algebra}

In this section, we recall the definition of W-algebras ($\hbar$-localized
DQ-algebras), and construct a W-algebra 
on $X$ by quantum Hamiltonian reduction. We introduce a quantized
symplectic coordinates of the W-algebra on $X$.
In the rest of the paper, we consider complex manifolds equipped with
the analytic topology. For a manifold $M$, we denote the sheaf of 
holomorphic functions on $M$ by $\calO_M$.

\subsection{Definition of W-algebras}
\label{sec:defin-w-algebr}

Let $\hbar$ be an indeterminant. Given $m \in \Z$, let 
$\W_{T^* \C^n} (m)$ be a sheaf of formal series 
$\sum_{k \ge -m} \hbar^k a_k$ ($a_k \in \calO_{T^* \C^n}$)
on the cotangent bundle $T^* \C^n$ of $\C^n$. We set
$\W_{T^* \C^n} = \bigcup_{m} \W_{T^* \C^n} (m)$. We define
a noncommutative $\C((\hbar))$-algebra structure on $\W_{T^* \C^n}$ by
\[
 f \circ g = \sum_{\alpha \in \Z_{\ge 0}^n} 
 \hbar^{|\alpha|} \frac{1}{\alpha!} \partial^{\alpha}_{\xi} f
 \cdot \partial^{\alpha}_{x} g
\]
where, for a multi-power $\alpha = (\alpha_1, \dots, \alpha_n) \in 
\Z^n_{\ge 0}$, $\alpha! = \alpha_1 ! \cdots \alpha_n !$ and
$|\alpha| = \alpha_1 + \dots + \alpha_n$.

Let $X$ be a complex symplectic manifold with symplectic form $\omega$.
A W-algebra on $X$ is a sheaf of $\C((\hbar))$-algebras $\W$ such that
for any point $x \in X$, there is an open neighbourhood $U$ of $x$,
a symplectic map $\phi: U \rightarrow T^* \C^n$, and a $\C((\hbar))$-algebra
isomorphism $\psi: \W \vert_U \isoto \phi^{-1} \W_{T^* \C^n}$.

The following fundamental properties of a W-algebra $\W$ are listed in
\cite{KR}.

\begin{enumerate}
\item 
The algebra $\W$ is a coherent and noetherian algebra.
\item
$\W$ contains a canonical $\C[[\hbar]]$-subalgebra $\W(0)$
which is locally isomorphic to $\W_{T^*\C^n}(0)$ (via the maps $\psi$).
We set $\W(m)=\hbar^{-m}\W(0)$.
\item
We have a canonical $\C$-algebra isomorphism
$\W(0)/\W(-1)\isoto \calO_X$ (coming from the canonical isomorphism via the maps
$\psi$).
The corresponding morphism
$\sigma_m: \W(m)\to \hbar^{-m}\calO_X$ is called the {\em symbol map}.
\item We have
\[
\sigma_{0}(\hbar^{-1}[f,g])=\{\sigma_0(f),\sigma_0(g)\}
\]
for any $f$, $g\in\W(0)$. Here $\{\scbul,\scbul\}$
is the Poisson bracket.
\item
The canonical map 
$\W(0)\to \prolim[{m\to\infty}]\W(0)/\W(-m)$
is an isomorphism.
\item
A section $a$ of $\W(0)$ is invertible in $\W(0)$ if and only if
$\sigma_0(a)$ is invertible in $\calO_X$.
\item
Given $\phi$, a $\C((\hbar))$-algebra automorphism of $\W$,
we can find locally an invertible section $a$ of $\W(0)$
such that $\phi=\Ad(a)$.
Moreover $a$ is unique up to a scalar multiple.
In other words, we have canonical isomorphisms
\[
\xymatrix{
\W(0)^\times/\C[[\hbar]]^\times\ar[r]^\sim_{\Ad}\ar[d]_\sim
&\Aut(\W(0))\ar[d]^\sim\\
\W^\times/\C((\hbar))^\times\ar[r]^\sim_{\Ad}&\Aut(\W).
} 
\]
\item
Let $v$ be a $\C((\hbar))$-linear filtration-preserving
derivation of $\W$.
Then there exists locally a  section $a$ of $\W(1)$
such that $v=\ad(a)$.
Moreover $a$ is unique up to a scalar.
In other words, we have an isomorphism
\[
\W(1)/\hbar^{-1}\C[[\hbar]]\isoto[\ad]\Der_{\mathrm{filtered}}(\W).
\]
\item
If $\W$ is a W-algebra, then its opposite ring
$\W^\opp$ is a W-algebra on $X^\opp$ where $X^\opp$ is
a symplectic manifold with symplectic form $- \omega$.
\end{enumerate}

A tuple $(f_1, \dots, f_n; g_1, \dots, g_n)$ of elements 
$f_i$, $g_i \in \W(0)$ is called quantized symplectic coordinates
of $\W$ if they satisfy $[f_i, f_j] = [g_i, g_j] = 0$ and
$[g_i, f_j] = \hbar \delta_{ij}$. 

For a $\W$-module $\scM$, a $\W(0)$-lattice of $\scM$ 
is a coherent $\W(0)$-submodule $\scM(0)$ such that
the canonical homomorphism $\W \otimes_{\W(0)} \scM(0) \rightarrow \scM$ 
is an isomorphism. We say that a $\W$-module $\scM$ is good
if, for any relatively compact open subset $U$ of $X$, there is a
coherent $\W(0)\vert_U$-lattice  of $\scM\vert_U$. We denote the
category of $\W$-modules by $\Mod(\W)$ and the 
full subcategory of good $\W$-modules by $\Mod^{good}(\W)$. Then,
$\Mod^{good}(\W)$ is an abelian subcategory of $\Mod(\W)$.

Next, we review the notion of F-actions.

Let $X$ be a symplectic manifold with the action of $\bbG_m$:
$\C^* \ni t \mapsto T_t \in \Aut(X)$. We assume there exists
a positive integer $m \in \Z_{> 0}$ such that 
$T^*_t \omega = t^m \omega$ for all $t \in \C^*$.

An F-action with exponent $m$ on $\W$ is an action of 
$\bbG_m$ on the
$\C$-algebra $\W$,
$\scF_t : T^{-1}_t \W \isoto \W$ for $t \in \C^*$ such that
$\scF_t(\hbar) = t^m \hbar$ and $\scF_t(f)$ depends holomorphically
on $t$ for any $f \in \W$.

An F-action with exponent $m$ on $\W$ extends to an F-action
with exponent $1$ on 
$\W[\hbar^{1/m}] = \C((\hbar^{1/m})) \otimes_{\C((\hbar))} \W$ given by
$\scF_t(\hbar^{1/m}) = t^1 \hbar^{1/m}$. 

\begin{definition}
 A $\W[\hbar^{1/m}]$-module with an F-action is a $\bbG_m$-equivariant
$\W[\hbar^{1/m}]$-module: i.e. there exist isomorphisms
$\scF_t: T^{-1}_t \scM \isoto \scM$ for $t \in \C^*$, and
we assume that
\begin{enumerate}
 \item $\scF_t(u)$ depends holomorphically on $t$ for any 
       $u \in \scM$;
 \item $\scF_t(f u) = \scF_t(f) \scF_t(u)$ for $f \in \W[\hbar^{1/m}]$
       and $u \in \scM$; and
 \item $\scF_t \circ \scF_{t'} = \scF_{t t'}$ for $t$, $t' \in \C^*$.
\end{enumerate}
\end{definition}

We denote by $\Mod_F(\W[\hbar^{1/m}])$ the category of 
$\W[\hbar^{1/m}]$-module with F-action, and by $\Mod_F^{good}(\W[\hbar^{1/m}])$
its full subcategory of good $\W[\hbar^{1/m}]$-modules with F-action.
These are $\C$-linear abelian categories. 

\subsection{Holonomic $\W$-modules}
\label{sec:holonomic-w-modules}

In this section, we review the notion of holonomic
$\W$-module introduced in \cite{KS-rh}.
The following proposition is due to \cite{KS}.

\begin{proposition}[\cite{KS}, Prop.~2.3.17]
\label{prop:5}
 For a good $\W$-module $\scM$, $\Supp \scM$ is involutive
 with respect to the Poisson bracket of $X$. In particular, 
 we have $\Dim \Supp \scM \ge \Dim X / 2$. 
\end{proposition}

\begin{definition}[\cite{KS-rh}]
 \label{def:8}
 We call a good $\W$-module $\scM$ holonomic
 if $\Supp \scM$ is a Lagrangian subvariety of $X$, i.e., 
 $\Dim \Supp \scM = \Dim X / 2$.
\end{definition}

\begin{proposition}
 \label{prop:6}
 The category of all holonomic $\W$-modules 
 $\Mod^{hol}(\W)$ is an abelian subcategory of 
 $\Mod^{good}(\W)$. 
\end{proposition}

The following lemma is obvious.

\begin{lemma}
 \label{lemma:4}
 Let $\scM$ be a good $\W$-module such that 
 $\Supp \scM$ is the disjoint union of subsets $Z_1$ and $Z_2$.
 Then, there exist two submodules $\scN_1$, $\scN_2$ of $\scM$,
 and we have $\scM = \scN_1 \oplus \scN_2$.
 \begin{proof}
  Define the submodules
  \[
   \scN_i = \{ m \in \scM \;\vert\; \Supp m \subset Z_i \}
  \]
  for $i=1$, $2$. Then the claim of lemma immediately follows.
 \end{proof}
\end{lemma}

In the present paper, we consider the case $\Dim X = 2$. 

Let $\bar{x}$, $\bar{\xi} \in \calO_{T^* \C^1}$ be coordinate functions
on $T^* \C^1$ defined by $\bar{x}((a, b)) = a$, $\bar{\xi}((a, b)) = b$
for $(a, b) \in T^* \C^1$. 
Let $x$, $\xi \in \W_{T^* \C^1}(0)$ be the standard quantized symplectic 
coordinates. That is we have $[\xi, x] = \hbar$ and $\sigma_0(x) = \bar{x}$,
$\sigma_0(\xi) = \bar{\xi}$. 

For $\lambda \in \C$, let $\scM_\lambda$ be a $\W_{T^* \C^1}$-module defined by
 \[
  \scM_\lambda = \W_{T^* \C^1} / \W_{T^* \C^1} (x \xi - \hbar \lambda).
 \]
Then, $\scM_\lambda$ is a holonomic $\W_{T^* \C^1}$-module supported
on $\{ \bar{x} \bar{\xi} = 0 \} \subset T^* \C^1$. Let $v_\lambda$ be the image
of the constant section $1 \in \W_{T^* \C^1}$ in $\scM_\lambda$.

\begin{lemma}
\label{lemma:15}
 For $m \in \Z$, we have the following isomorphism of 
 $\W_{T^* \C^1} \vert_{\{\bar{x} \ne 0\}}$-modules:
 \begin{align*}
  \scM_\lambda \vert_{\{\bar{x} \ne 0\}} &\longrightarrow
  \scM_{\lambda + m} \vert_{\{\bar{x} \ne 0 \}}, \\
  v_{\lambda} &\mapsto x^{-m} v_{\lambda+m}.
 \end{align*}
 Obviously, the inverse homomorphism is given by 
 $v_{\lambda+m} \mapsto x^m v_\lambda$.
\end{lemma}

A similar proposition holds globally on $T^* \C^1$. 
It is an analogue of a well-known fact on
regular holonomic $\calD_{\C^1}$-modules.

\begin{proposition}
\label{prop:11}
 For any $\lambda \ne -1$, we have an isomorphism of 
 $\W_{T^* \C^1}$-modules, $\scM_\lambda \simeq \scM_{\lambda+1}$.
 \begin{proof}
  Define homomorphisms of $\W_{T^* \C^1}$-modules
  \begin{equation*}
   \phi: \scM_\lambda \longrightarrow \scM_{\lambda+1}, \quad
   v_\lambda \mapsto \hbar^{-1} \xi v_{\lambda+1},
  \end{equation*}
  and
  \begin{equation*}
   \psi: \scM_{\lambda+1} \longrightarrow \scM_\lambda, \quad
   v_{\lambda+1} \mapsto \frac{1}{\lambda+1} x v_\lambda.
  \end{equation*}
  These homomorphisms are mutually inverse, i.e.
  \begin{equation*}
   \phi \circ \psi (v_{\lambda+1}) = \phi\left(\frac{1}{\lambda+1} x v_\lambda \right)
   = \frac{\hbar^{-1}}{\lambda+1} (x \circ \xi) v_{\lambda+1} = v_{\lambda+1},
  \end{equation*}
  and
  \begin{equation*}
   \psi \circ \phi (v_\lambda) = \psi(\hbar^{-1} \xi v_{\lambda})
   = \frac{\hbar^{-1}}{\lambda+1} (\xi \circ x) v_{\lambda} = v_{\lambda}.
  \end{equation*}
  Therefore $\scM_\lambda$ and $\scM_{\lambda+1}$ are isomorphic. 
 \end{proof}
\end{proposition}

The following proposition is essential for the microlocal analysis
of holonomic $\W_{T^* \C^1}$-modules. This is an analogue
of a consequence of the classification theorem of simple holonomic system 
(cf. \cite[Proposition~8.36]{Ka}). 


\begin{proposition}
 \label{prop:7}
 Set $Z_1 = \{x=0\}$ and $Z_2 = \{\xi=0\}$.
 Note that $\Supp \scM_\lambda = Z_1 \cup Z_2$. Then, 
 \begin{enumerate}
  \item For $\lambda \not\in \Z$, $\scM_\lambda$ is an irreducible
	$\W_{T^* \C^1}$-module.
  \item For $\lambda \in \Z_{\ge 0}$, there exist a 
	$\W_{T^* \C^1}$-submodule $\scN$ of $\scM_\lambda$ supported
	on $Z_1$,
	and $\Supp \scM_\lambda / \scN = Z_2$  on a neighborhood of $\{x = \xi = 0\}$.
  \item For $\lambda \in \Z_{< 0}$, there exist a 
	$\W_{T^* \C^1}$-submodule $\scN$ of $\scM_\lambda$ supported
	on $Z_2$,
	and $\Supp \scM_\lambda / \scN = Z_1$
	on a neighborhood of $\{x = \xi = 0\}$.
 \end{enumerate}
 \begin{proof}
  The proof of this proposition is similar to one of 
  \cite[Proposition 8.36]{Ka}. 




 \end{proof}
\end{proposition}

\subsection{W-algebra $\tildA_c$ on the quiver variety $X$}
\label{sec:w-algebra-on-quiver-var}

We denote the restriction of the canonical W-algebra $\W_{T^* \C^{l}}$
to $\widetilde{X} \subset T^* \C^{l}$ by $\W_{\widetilde{X}}$.
Let $(x_1, \dots, x_l; y_1, \dots, y_l)$ ($x_i$, $y_i \in \W_{T^* \C^l}$)
be the standard quantized symplectic coordinates: 
$[x_i, x_j] = [y_i, y_j] = 0$ and $[y_i, x_j] = \delta_{ij} \hbar$ for
all $i$, $j$.
The action of the reductive group $G$ on $\widetilde{X}$ induces
an action on $\W_{\widetilde{X}}$. We define
the following homomorphism $\mu_{\widetilde{X}}$ of Lie algebras
\begin{align*}
 \mu_{\widetilde{X}}: \g &\longrightarrow \W_{\widetilde{X}}(1), \\
 A_i &\mapsto \hbar^{-1} (x_{i+1} y_{i+1} - x_i y_i).
\end{align*}
We call $\mu_{\widetilde{X}}$ a quantum moment map with respect to
the action of $G$. Fix a parameter $c_0$, $\dots$, $c_{l-1} \in \C$ such
that $c_0 + \dots + c_{l-1} = 0$.
We define a $\W_{\widetilde{X}}$-module 
$\scL_{c}$ by
\[
 \scL_c = \W_{\widetilde{X}} \Bigm/ 
 \sum_{i=0}^{l-1} \W_{\widetilde{X}} (\mu_{\widetilde{X}}(A_i) + c_i)
 = \W_{\widetilde{X}} \Bigm/ 
 \sum_{i=0}^{l-1} \W_{\widetilde{X}} (x_{i+1} y_{i+1} - x_i y_i + \hbar c_i).
\]
The $\W_{\widetilde{X}}$-module $\scL_c$ is a good 
$\W_{\widetilde{X}}$-module with a $\W_{\widetilde{X}}(0)$-lattice
\[
\scL_c(0) \defeq \W_{\widetilde{X}}(0) \Bigm/ \sum_{i=0}^{l-1}
\W_{\widetilde{X}}(0) (x_{i+1} y_{i+1} - x_i y_i + \hbar c_i).
\]
Define a sheaf of algebras on $X$,
\[
 \A_c = \left(p_{*} \lEnd_{\W_{\widetilde{X}}}(\scL_c)^{G}\right)^{\opp}
\]
where $p: \mu^{-1}(0) \rightarrow X$ is the projection.
By \cite{KR}, $\A_c$ is a W-algebra on $X$. Set
\[
 \A_c(0) = \left(p_{*} \lEnd_{\W_{\widetilde{X}}(0)}(\scL_c(0))^{G}\right)^{\opp}.
\]
Then, $\A_c(0)$ is a canonical $\C[[\hbar]]$-subalgebra of $\A_c$. 

Let us give an F-action on $\W_{\widetilde{X}}$ by $\scF_t(x_i) = t x_i$,
$\scF_t(y_i) = t y_i$, and $\scF_t(\hbar) = t^2 \hbar$ for 
$t \in \C^*$. The corresponding $\bbG_m$-action on $\widetilde{X}$
is given by $T_t ((a_i, b_i)_{i=1,\dots,l}) = (t a_i, t b_i)_{i=1,\dots,l}$.
This action induces a $\bbG_m$-action on the quiver variety $X$ which
coincides with the action induced from the toric $\bbT^2$-action on $X$.

The F-actions on $\W_{\widetilde{X}}$ induce an F-action
with exponent $2$ on $\A_c$.
We set $\tildA_c = \A_c[\hbar^{1/2}]$ and
$\tildA_c(0) = \A_c(0)[\hbar^{1/2}]$.



In \cite{BK}, we have the following W-affinity of the algebra 
$\tildA_c$.

\begin{theorem}[\cite{BK}]
 \label{thm:1}
 Let 
 $A_c = \left(
\gEnd_{\Mod_{F}^{good}(\tildA_c)}(\tildA_c)\right)^{\opp}$.

 Assume that $c_i + c_{i+1} + \dots + c_{j-1} \in \Z_{\ge 0}$
 implies $\theta_i + \theta_{i+1} + \dots + \theta_{j-1} < 0$,
i.e. $i \rhd j$.
Then
 we have the following equivalence of categories:
 \begin{align*}
 \Mod_{F}^{good}(\tildA_c) &\simeq 
 A_c\mmod, \\
 \scM &\mapsto \gHom_{\Mod_{F}^{good}(\tildA_c)}(\tildA_c, \scM).
 \end{align*}
 Its quasi-inverse functor is given by 
 $M \mapsto \tildA_c \otimes_{A_c} M$.
\end{theorem}

\begin{remark}
\label{rmk:1}
 As we see in \refsec{sec:RCA}, 
 the algebra $A_c$ is isomorphic to the spherical subalgebra of 
 the rational Cherednik algebra of type $\Z / l \Z$.
\end{remark}

In the rest of this paper, we assume the assumptions of \refthm{thm:1}.

Let $\One_c$ be the image of the constant section $1 \in \W_{\widetilde{X}}$ 
in $\scL_c$. For a $G$-invariant section $f \in \W_{\widetilde{X}}$,
a $G$-invariant endomorphism of $\scL_c$ is uniquely defined by the right multiplication
$g \One_c \mapsto g f \One_c \in \lEnd_{\W_{\widetilde{X}}}(\scL_c)$
for $g \in \W_{\widetilde{X}}$. By abuse of notation,
we denote the image of the above $G$-invariant endomorphism in
$\tildA_c$ by the same symbol $f$.

Consider global sections $x_1 \cdots x_l$, $y_1 \cdots y_l$,
$x_i y_i$ of $\tildA_c$ for $i=1$, $\dots$, $l$. Although these
sections  are not F-invariant, the global sections 
$\hbar^{-l/2} x_1 \cdots x_l$, $\hbar^{-l/2} y_1 \cdots y_l$,
$\hbar^{-1} x_i y_i$ are elements of 
$A_c = \left(\gEnd_{\Mod_F^{good}(\tildA_c)}(\tildA_c)\right)^{\opp}$.
In $\tildA_c$, we have relations
\[
 x_{i+1} y_{i+1} - x_i y_i + \hbar c_i = 0
\]
for $i=1$, $\dots$, $l$.

Next, we consider the local structure of the W-algebra 
$\tildA_c$ on the affine open subset $X_i$ for
$i=1$, $\dots$, $l$. 
Set $\tildA_{c,i} = \tildA_c \vert_{X_i}$.
Recall the indices $\ii_1$, $\dots$, $\ii_l$ in \refeq{eq:2}.
We define local sections of $\tildA_c$
\[
 f_i = (x_{\ii_1} \cdots x_{\ii_i}) \circ (y_{\ii_{i+1}} \cdots y_{\ii_l})^{-1}, \quad
 g_{i} = (y_{\ii_{i}} \cdots y_{\ii_l}) \circ (x_{\ii_1} \cdots x_{\ii_{i-1}})^{-1},
\]
on $X_i$. 
We have $f_i \circ g_i = x_{\ii_i} y_{\ii_i}$ and 
$g_i \circ f_i = x_{\ii_i} y_{\ii_i} + \hbar$. 
Thus, 
for $i=1$, $\dots$, $l$, 
$\tildA_{c,i}$ is isomorphic to
$\W_{T^* \C^1}$ by
$x \mapsto f_i$, $\xi \mapsto g_i$.
That is, $(f_i; g_i)$ is a quantized symplectic coordinate
of $\tildA_{c,i}$.
We have 
\begin{equation}
 \label{eq:26}
  y_1 \cdots y_l = g_i \circ (x_{\ii_1} y_{\ii_1}) \circ \cdots 
  \circ (x_{\ii_{i-1}} y_{\ii_{i-1}}) \quad \text{on } X_i
\end{equation}
We also have
$g_{i+1} \circ f_i = f_i \circ g_{i+1} = 1$ in $\tildA_c \vert_{X_i \cap X_{i+1}}$. 
Sometimes, we denote the section $g_{i+1} \vert_{X_i \cap X_{i+1}}$ by $f_{i}^{-1}$. 

For $i=1$, $\dots$, $l-1$, we set 
\[
\tilde{c}_i = c_{\ii_i} + c_{\ii_{i} + 1} + \dots + c_{\ii_{i+1} - 1}.
\]
Under the assumption of \refthm{thm:1}, we have 
$\tilde{c}_i + \tilde{c}_{i+1} + \dots + \tilde{c}_{j-1} \not\in \Z_{\le 0}$
for $1 \le i < j \le l$.
Then, we have
\begin{equation}
\label{eq:4}
 x_{\ii_{i+1}} y_{\ii_{i+1}} - x_{\ii_i} y_{\ii_i} + \hbar \tilde{c}_i = 0.
\end{equation}




\section{The rational Cherednik algebra and category $\calO$}
\label{sec:RCA}

In this section, we review the definition and fundamental 
facts of the rational Cherednik algebra of type $\Z / l \Z$
and the category $\calO$. 

Let $\Z / l\Z = \langle \gamma \rangle$ be a 
cyclic group with an action on $\C$ given by 
$\gamma \mapsto \zeta = \exp(2\pi\sqrt{-1}/l)$. Let 
$D(\C^*)$ be the algebra of algebraic differential operators on $\C^*$.
Let $z$ be a standard coordinate function on $\C$. Then, we have
$\C = \Spec[z]$ and $\C^* = \Spec\C[z, z^{-1}]$. The algebra
$D(\C^*)$ is generated by $z^{\pm}$ and $d / dz$.

The action of $\Z/l\Z$ on $\C$ induces an action of $D(\C^*)$
given by $\gamma(z) = \zeta^{-1} z$, $\gamma(d/dz) = \zeta d/dz$.
We denote the smash product of $D(\C^*)$ and $\Z/l\Z$ by
$D(\C^*) \# \Z/l\Z$. 

For a parameter $\kappa = (\kappa_1, \dots, \kappa_{l-1}) \in \C^{l-1}$,
we define a Dunkl operator $\partial_\kappa$ by
\[
 \partial_{\kappa} = \frac{d}{dz} + \frac{l}{z} \sum_{i=0}^{l-1} \kappa_i 
 \bfe_i
\]
where we regard $\kappa_0 = 0$ and let 
$\bfe_i = (1/l) \sum_{j=0}^{l-1} \zeta^{ij} \gamma^j$ be an
idempotent of $\C (\Z/l\Z)$ for $i=0$, $1$, $\dots$, $l-1$.

\begin{definition}[\cite{EG}]
 \label{def:10}
 \begin{enumerate}
  \item The rational Cherednik algebra $H_\kappa = H_\kappa(\Z/l\Z)$ is 
	the subalgebra of $D(\C^*) \# \Z/l\Z$ generated by 
	$z$, $\partial_\kappa$ and $\gamma$.
  \item The spherical subalgebra of $H_\kappa$ is the algebra
	$\bfe_0 H_\kappa \bfe_0$.
 \end{enumerate}
\end{definition}

The following proposition is analogue of the triangular decomposition
of semisimple Lie algebras.

\begin{proposition}[\cite{EG}]
\label{prop:12}
 We have the following isomorphisms as $\C$-linear spaces
 \[
  H_\kappa \simeq \C[z] \otimes_{\C} \C (\Z/l\Z) \otimes_{\C} \C[\partial_\kappa], \quad \text{and} \quad
 \bfe_0 H_\kappa \bfe_0 \simeq \C[z, \partial_\kappa]^{\Z/l\Z}
 = \C[z^l, z \partial_{\kappa}, \partial_{\kappa}^l].
 \]
\end{proposition}

By \cite{BK} together with \cite{Ho}, we have the following isomorphism
of algebras (\refremark{rmk:1})
\begin{align}
 \label{eq:16}
 \bfe_0 H_\kappa \bfe_0 &\longrightarrow A_c, \\
 \bfe_0 z^l \bfe_0 &\mapsto \hbar^{-l/2} x_1 \cdots x_l, \nonumber\\
 \bfe_0 \partial_\kappa^l \bfe_0 &\mapsto \hbar^{-l/2} y_1 \cdots y_l, \nonumber\\
 \bfe_0 z \partial_\kappa \bfe_0 &\mapsto \hbar^{-1} x_1 y_1. \nonumber
\end{align}
where 
\begin{equation}
\label{eq:17} 
c = c(\kappa) = (c_i)_{i=0, 1, \dots, l-1}, \qquad
c_i = \kappa_{i} - \kappa_{i+1} - 1/l + \delta_{i, 0}.
\end{equation}

\begin{remark}
 In \cite{K}, we use parameters $\lambda_i = - c_i$ for quantum 
 Hamiltonian reduction.
\end{remark}


\begin{lemma}[cf. \cite{K}, Proposition 4.4] 
\label{lemma:13}
 Assume $c_i + c_{i+1} + \dots + c_{j-1} \ne 0$ for
 $0 < i < j \le l$.
Then, the rational Cherednik algebra
 $H_\kappa$ is Morita equivalent to its spherical subalgebra
 $\bfe_0 H_\kappa \bfe_0 \simeq A_c$, i.e. we have an equivalence
 of categories
 \begin{align*}
  H_\kappa\mmod &\longrightarrow (\bfe_0 H_\kappa \bfe_0)\mmod, \\
  M &\mapsto \bfe_0 M.
 \end{align*}
\end{lemma}

In the present paper, we assume the assumption of \reflemma{lemma:13}
holds.

The category $\calO(H_\kappa)$ is the subcategory of $H_\kappa\mmod$
such that the Dunkl operator $\partial_\kappa$ acts locally nilpotently
on a module $M \in \calO(H_\kappa)$. 

Consider an irreducible $\C (\Z/l\Z)$-module $\C \bfe_i$ for $i=0$, $\dots$, $l$.
We regard $\C \bfe_i$ as a $\C[\partial_\kappa] \# \Z/l\Z$-module by
$\partial_\kappa \bfe_i = 0$. We define an $H_\kappa$-module
\[
 \bfDelta_\kappa(i) = H_\kappa \otimes_{\C[\partial_\kappa] \# \Z/l\Z}
 \C \bfe_i
\]
called a standard module. By \refprop{prop:12}, we have
\begin{equation}
\label{eq:20} 
 \bfDelta_\kappa(i) = \C[z] \bfe_i
\end{equation}
as a $\C$-linear space. 

By the equivalence of \reflemma{lemma:13}, we have a subcategory
$\calO(A_c)$ of $A_c\mmod$ which is equivalent to the category
$\calO(H_\kappa)$. We call $\calO(A_c)$ the category $\calO$ of $A_c$.
For $i=1$, $\dots$, $l$, an $A_c$-module 
\[
 \Delta_c(i) = \bfe_0 \bfDelta_\kappa(i)
\]
where $c$ is given by \refeq{eq:17} and we regard $\bfe_l = \bfe_0$.
The module $\Delta_c(i)$ is a standard module for $\calO(A_c)$.

The following proposition is fundamental and well-known facts about
the category $\calO(A_c)$ and modules $\Delta_c(i)$ ($i=1$, $\dots$,
$l$).

\begin{proposition}[\cite{GGOR}]
 We have the following fundamental facts about the standard modules
 $\Delta_c(i)$:
\label{prop:3}
 \begin{enumerate}
  \item For $i=1$, $\dots$, $l$, the standard module $\Delta_c(i)$
	has a unique irreducible quotient $L_c(i)$.
  \item The irreducible modules $L_c(i)$ ($i=1$, $\dots$, $l$) are
	mutually non-isomorphic.
  \item Any simple object in the category $\calO(A_c)$ is isomorphic
	to $L_c(i)$ for some $i=1$, $\dots$, $l$.
 \end{enumerate}
\end{proposition}

\begin{remark}
 Originally \cite{GGOR} considered the category $\calO(H_\kappa)$ of the
 rational Cherednik algebra $H_\kappa$, 
 not one of its spherical subalgebra $\bfe_0 H_\kappa \bfe_0 \simeq A_c$.

\end{remark}

By \refeq{eq:20}, we have 
\begin{equation}
 \label{eq:22}
 \Delta_c(i) = \bfe_0 \bfDelta_\kappa(i) = \bfe_0 \C[z^l] z^{l-i} \bfe_i
 = \C[\hbar^{-l/2} x_1 \cdots x_l] e_i
\end{equation}
as a $\C$-linear space where we denote $\bfe_0 z^{l-i} \bfe_i$ by $e_i$.

In $A_c$, we have
\begin{align*}
 [\hbar^{-1} x_{\ii_1} y_{\ii_1}, \hbar^{-l/2} x_1 \cdots x_l] &= 
 \hbar^{-l/2} x_1 \cdots x_l, \\
 [\hbar^{-1} x_{\ii_1} y_{\ii_1}, \hbar^{-l/2} y_1 \cdots y_l] &= - 
 \hbar^{-l/2} y_1 \cdots y_l.
\end{align*}
For $i=1$, $\dots$, $l$,
the operator $\hbar^{-1} x_{\ii_1} y_{\ii_1}$ acts
semisimply on the standard module $\Delta_c(\ii_i)$, i.e. $\Delta_c(\ii_i)$
is a direct sum of eigenspaces with respect to the action of 
$\hbar^{-1} x_{\ii_1} y_{\ii_1}$. In fact, by direct calculation we have
\[
 \Delta_c(\ii_i) = \bigoplus_{m \in \Z_{\ge 0}} (\hbar^{-l/2} x_1 \cdots x_l)^m e_{\ii_i},
\]
and 
\begin{equation}
\label{eq:23} 
 (\hbar^{-1} x_{\ii_1} y_{\ii_1}) \circ (\hbar^{-l/2} x_1 \cdots x_l)^m e_{\ii_i} = 
 (m + \tilde{c}_{\ii_1} + \tilde{c}_{\ii_2} + \dots + \tilde{c}_{\ii_{i-1}})
 (\hbar^{-l/2} x_1 \cdots x_l)^m e_{\ii_i}.
\end{equation}

\begin{lemma}
 \label{lemma:14}
 We have 
\begin{equation*}
  \Delta_c(\ii_i) = A_c / (A_c (\hbar^{-1} x_{\ii_i} y_{\ii_i}) + A_c (\hbar^{-l/2} y_1 \cdots y_l))
\end{equation*}
for $i=1$, $\dots$, $l$.
 \begin{proof}
The standard module $\Delta_c(\ii_i)$ is cyclic with cyclic vector $e_{\ii_i}$.
By \refeq{eq:23}, we have $\hbar^{-1} x_{\ii_i} y_{\ii_i} e_{\ii_i} = 0$.
Thus, we have the following surjective homomorphism of $A_c$-modules
\begin{align*}
 A_c / (A_c (\hbar^{-1} x_{\ii_i} y_{\ii_i}) + A_c (\hbar^{-l/2} y_1 \cdots y_l)) 
 &\twoheadrightarrow \Delta_c(\ii_i), \\
 f &\mapsto f e_{\ii_i}.
\end{align*}
By \refprop{prop:12} and \refeq{eq:22}, this homomorphism is an isomorphism.
 \end{proof}
\end{lemma}





\section{Microlocal construction of modules}
\label{sec:micr-constr-modul}

\subsection{Construction of the standard modules}
\label{sec:constr-stand-mod}

In this section, we introduce the $\tildA_c$-module
$\calM_{c}^{\Delta}(\ii_i)$ supported
on a Lagrangian subvariety
 $D_i \cup D_{i+1} \cup \dots \cup D_l$. Moreover, we show that
the $\calM_c^{\Delta}(\ii_i)$ is 
a counterpart of the standard module $\Delta_c(\ii_i)$
of $A_c$ through the 
equivalence of \refthm{thm:1}.

\begin{definition}
 \label{def:9}
 For $1 \le i < i' \le l$ and a parameter
 $\lambda = (\lambda_j)_{j=i+1, \dots, i'} \in \C^{i'-i}$,
 we call $\lambda$ admissible when
 it satisfies $\lambda_j - \lambda_{j+1} - \tilde{c}_j \in \Z$
 for $j=i, \dots, i'-1$ where we regard $\lambda_i = 0$.
\end{definition}

\begin{definition}
\label{def:6}
 For $i=1$, $\dots$, $l$, we 
 take an admissible parameter $\lambda = (\lambda_{j})_{j=i+1, \dots, l}$.
 We define an $\tildA_c$-module 
 $\calM_{c,\lambda}(\ii_i)$ by glueing
 local sheaves as follows:
  \begin{align*}
  \calM_{c, \lambda}(\ii_i) \vert_{X_i} &= 
  \tildA_{c, i} / \tildA_{c, i} g_i, \\
  {\calM_{c, \lambda}(\ii_i)} \vert_{X_j} &= 
  \tildA_{c,j} / \tildA_{c,j} (f_j \circ g_j - \hbar \lambda_j) \\
   &= \tildA_{c,j} / \tildA_{c,j} (x_{\ii_j} y_{\ii_j} - \hbar \lambda_j) 
  \qquad (\text{for } j = i+1, \dots, l), \\
  {\calM_{c, \lambda}(\ii_i)} \vert_{X_j} &= 0 
  \qquad (\text{for } j = 1, \dots, i-1).
 \end{align*}
  Its glueing is given by 
\begin{equation}
\label{eq:6}
  u_j = f_{j}^{\lambda_j - \lambda_{j+1} - \tilde{c}_j} 
 u_{j+1} \quad \text{on } X_j \cap X_{j+1}
\end{equation}
 where $u_j$ is the image of the constant section $1 \in \A_{c,j}$ in
 $\calM_{c,\lambda}(\ii_i) \vert_{X_j}$ for $j=i$, $\dots$, $l$.
\end{definition}

Note that we have 
\begin{equation}
 \label{eq:12}
 \calM_{c, \lambda}(\ii_i) \vert_{X_j} \simeq \scM_{\lambda_j}, \qquad
 u_j \mapsto v_{\lambda_j} 
\end{equation}
under the isomorphism $\tildA_{c,j} \simeq \W_{T^* \C^1}$.

\begin{lemma}
\label{lemma:9}
 The module $\calM_{c,\lambda}(\ii_i)$ is a well-defined 
 good $\tildA_c$-module supported on the Lagrangian subvariety
 $D_i \cup D_{i+1} \cup \dots \cup D_{l}$.
 \begin{proof}
  Set $\scN_1 = \calM_{c, \lambda}(\ii_i) \vert_{X_j}$ and
  $\scN_2 = \calM_{c, \lambda}(\ii_i \vert_{X_{j+1}})$. By
  \refeq{eq:4}, we have
  \[
   f_{j+1} \circ g_{j+1} - f_j \circ g_j + \hbar \tilde{c}_j = 0
  \]
  on $X_j \cap X_{j+1}$. Thus, we have
\begin{align*}
   \scN_2 \vert_{X_j \cap X_{j+1}} 
 &= \tildA_c \vert_{X_j \cap X_{j+1}} \bigm/ 
 \tildA_c \vert_{X_j \cap X_{j+1}} (f_{j+1} \circ g_{j+1} - \hbar \lambda_{j+1}), \\
 &= \tildA_c \vert_{X_j \cap X_{j+1}} \bigm/ 
 \tildA_c \vert_{X_j \cap X_{j+1}} (f_{j} \circ g_{j} - \hbar (\lambda_{j+1} + \tilde{c}_j)).
\end{align*}
  Since $\lambda$ is admissible, 
  by \reflemma{lemma:15}, $\scN_1 \vert_{X_j \cap X_{j+1}}$ is
  isomorphic to $\scN_2 \vert_{X_j \cap X_{j+1}}$ by 
  $u_j \mapsto f_j^{\lambda_j - \lambda_{j+1} - \tilde{c}_j} u_{j+1}$,
  and consequently $\calM_{c, \lambda}(\ii_i)$ is well-defined.

  There exists an $\tildA_c(0)$-lattice of $\calM_{c, \lambda}(\ii_i)$ 
  given by $\calM_{c,\lambda}(0) \vert_{X_j} = \tildA_c(0) u_j$ for
  $j=i$, $\dots$, $l$.
  Thus, $\calM_{c, \lambda}(\ii_i)$ is a good $\tildA_c$-module.
 \end{proof}
\end{lemma}

\begin{lemma}
\label{lemma:10}
 Fix $i = 1$, $\dots$, $l$. Take any admissible parameter
 $\lambda = (\lambda_j)_{j=i+1, \dots, l} \in \C^{l-i}$ such that
 $\lambda \in (\C \backslash \Z_{< 0})^{l-i}$ 
 (resp. $\lambda \in (\C \backslash \Z_{\ge 0})^{l-i}$) and 
 $\lambda' \in \lambda + (\Z_{\ge 0})^{l-i}$
 (resp. $\lambda' \in \lambda + (\Z_{\le 0})^{l-i}$). 
 Then, we have an isomorphism of $\tildA_{c}$-modules
 $\calM_{c, \lambda}(\ii_i) \simeq \calM_{c,\lambda'}(\ii_i)$.
\begin{proof}
 We will prove the case where $\lambda \in (\C \backslash \Z_{< 0})^{l-i}$ and
 $\lambda' \in \lambda + (\Z_{\ge 0})^{l-i}$. 
  It is enough to show that the claim of the lemma holds when
  there exists $j = i+1$, $\dots$, $l$ such that
  $\lambda'_j + 1 = \lambda_j$ 
  and $\lambda'_k = \lambda_k$ for $k \ne j$. 
 
 By \refeq{eq:12}, we have
 \[
  \calM_{c,\lambda}(\ii_i)\vert_{X_j} \simeq \scM_{\lambda_j}, \qquad
 \calM_{c, \lambda'}(\ii_i)\vert_{X_j} \simeq \scM_{\lambda_j + 1}.
 \]
 Thus, there exists an isomorphism of $\tildA_{c,j}$-modules
 $\calM_{c,\lambda}(\ii_i)\vert_{X_j} \simeq \calM_{c,\lambda'}(\ii_i)\vert_{X_j}$ 
 by \refprop{prop:11}. For $k \ne j$, we have a trivial isomorphism
 of $\tildA_{c,k}$-modules 
 $\calM_{c,\lambda}(\ii_i)\vert_{X_k} \simeq \calM_{c,\lambda'}(\ii_i)\vert_{X_k}$.
 These isomorphisms induce an isomorphism of $\tildA_c$-module
 $\calM_{c, \lambda}(\ii_i) \simeq \calM_{c, \lambda'}(\ii_i)$.

 The case where $\lambda \in (\C \backslash \Z_{\ge 0})^{l-i}$ and
 $\lambda' \in \lambda + (\Z_{\le 0})^{l-i}$ is proved similarly. 
\end{proof}
\end{lemma}

Then, we define the $\tildA_c$-modules $\calM_c^{\Delta}(\ii_i)$, 
$\calM_c^{\nabla}(\ii_i)$.

\begin{definition}
 \label{def:7}
 For an admissible parameter
 $\lambda \in (\C \backslash \Z_{\ge 0})^{l-i}$, we denote 
  \begin{equation*}
   \calM_c^{\Delta}(\ii_i) =
   \calM_{c,\lambda}(\ii_i).
 \end{equation*}
\end{definition}

\begin{remark}
 For an admissible parameter
 $\lambda \in (\C \backslash \Z_{< 0})^{l-i}$, we denote
 \begin{equation*}
   \calM_c^{\nabla}(\ii_i) =
   \calM_{c,\lambda}(\ii_i).
 \end{equation*}
 The module $\calM_c^{\nabla}(\ii_i)$ is an $\tildA_c$-module
 (conjecturally) corresponding to a costandard module of
 $A_c$.
\end{remark}

In the rest of this section, we show that the $\tildA_c$-module
$\calM_c^{\Delta}(\ii_i)$ corresponds to the standard module
$\Delta_c(\ii_i)$ via the equivalence of categories \refthm{thm:1},
i.e. we have $\gHom_{\Mod_F^{good}(\tildA_c)}(\tildA_c, \calM_c^{\Delta}(\ii_i)) \simeq \Delta_c(\ii_i)$.


\begin{theorem}
 \label{thm:2}
 We have an isomorphism of $\tildA_c$-modules
 $\calM_{c}^{\Delta}(\ii_i) \simeq \tildA_c \otimes_{A_c} \Delta_c(\ii_i)$.
 In other words, we have an isomorphism of $A_c$-modules,
 \[
\gHom_{\Mod_F^{good}(\tildA_c)}(\tildA_c, \calM_c^{\Delta}(\ii_i)) \simeq \Delta_c(\ii_i).
 \]
  \begin{proof}
  By \reflemma{lemma:14}, we have
\begin{equation}
\label{eq:27} 
   \tildA_c \otimes_{A_c} \Delta_c(\ii_i) 
  \simeq \tildA_c / (\tildA_c x_{\ii_i} y_{\ii_i} + \tildA_c y_1 \cdots y_l).
\end{equation}

  For $j = 1$, $\dots$, $i$, we have on $X_j$
  \begin{align*}
   y_1 \cdots y_l &= g_j \circ x_{\ii_1} y_{\ii_1} \circ \dots \circ x_{\ii_{j-1}} y_{\ii_{j-1}} \\
   &= g_j \circ \prod_{k=1}^{j-1} (f_j \circ g_j + \hbar (\tilde{c}_k + \cdots +\tilde{c}_{j-1})).
  \end{align*}
  Since, $\hbar (\tilde{c}_k + \dots + \tilde{c}_{i-1})$ is nonzero 
  in the field $\C((\hbar))$ for $k=1$, $\dots$, $i-1$,
  the isomorphism \refeq{eq:27} is reduced on $X_j$ as follows
  \begin{align*}
   \lefteqn{\tildA_{c,j} \otimes_{A_c} \Delta_c(\ii_i)} & \\
   &\simeq \tildA_{c,j} \Bigm/ \bigl\{ \tildA_{c,j} (f_j \circ g_j - \hbar (\tilde{c}_j + \dots \tilde{c}_{i-1})), \\
   & \hspace*{4cm} \textstyle
   + \tildA_{c,j} \;g_j \circ \prod_{k=1}^{j-1} (f_j \circ g_j - \hbar (\tilde{c}_k + \dots + \tilde{c}_{j-1})) \bigr\}, \\
   &= \tildA_{c,j} \Bigm/ \bigl\{ \tildA_{c,j} (f_j \circ g_j - \hbar (\tilde{c}_j + \dots + \tilde{c}_{i-1})) + \tildA_{c,j} g_j \bigr\}, \\
   &= 
\begin{cases}
 \tildA_{c,i} / \tildA_{c,i} g_i & \text{for } j=i, \\
 0 & \text{for } j=1, \dots, i-1.
\end{cases}
  \end{align*}

For $j=i+1$, $\dots$, $l$, we have on $X_j$,
\begin{equation*}
 y_1 \cdots y_l = g_j \circ (x_{\ii_1} y_{\ii_1}) \circ \dots \circ
 (x_{\ii_i} y_{\ii_i}) \circ \dots \circ (x_{\ii_{j-1}} y_{\ii_{j-1}}) 
\end{equation*}
  by \refeq{eq:26}. Since $x_{\ii_k} y_{\ii_k}$ and $x_{\ii_i} y_{\ii_i}$
  commute with each other, on $X_j$ we have
  \begin{align*}
   \tildA_{c,j} \otimes_{A_c} \Delta_{c}(\ii_i) &\simeq
   \tildA_{c,j} \Bigm/ \Bigl( \tildA_{c,j} x_{\ii_i} y_{\ii_i} +
   \tildA_{c,j} g_j \circ \bigl(\textstyle \prod_{k \ne i} x_{\ii_k} y_{\ii_k}\bigr) \circ
   x_{\ii_i} y_{\ii_i}
   \Bigr), \\
   &= \tildA_{c,j} / \tildA_{c,j} (x_{\ii_j} y_{\ii_j} 
   + \hbar (\tilde{c}_i + \dots \tilde{c}_{j-1})) \qquad 
   \text{for } j=i+1, \dots, l.
  \end{align*}
  Note that $\lambda = (\lambda_j)_{j=i+1, \dots, l} \in \C^{l-i}$ where
  $\lambda_j = - (\tilde{c}_i + \dots + \tilde{c}_{j-1})$ is admissible
  and $\lambda \in (\C \backslash \Z_{\ge 0})^{l-i}$. Thus, 
  the $\tildA_c$-module $\tildA_c \otimes_{A_c} \Delta_c(\ii_i)$ is
  clearly isomorphic to $\calM_c^{\Delta}(\ii_i)$.
 \end{proof}
\end{theorem}

\subsection{Construction of irreducible modules of $\tildA_c$}
\label{sec:constr-irred-modul}

In this subsection, we construct modules $\calL_c(i)$ of the
W-algebra $\tildA_c$ for $i=1$, $\dots$, $l$, and
show that they are irreducible modules. Under the equivalence of
\refthm{thm:1}, 
$\gHom_{\Mod_{F}^{good}(\tildA_c)}(\tildA_c, \calL_c(i))$
is isomorphic to the irreducible module $L_c(i)$ of $A_c$ defined
in \refsec{sec:RCA}.

 Fix $i=1$, $\dots$, $l$. We denote by $\epsilon(i) = i+1$, $\dots$, $l+1$ 
 a unique index such that
 $\tilde{c}_i + \tilde{c}_{i+1} + \dots + \tilde{c}_{\epsilon(i)-1} \in \Z$ and
 $\tilde{c}_i + \tilde{c}_{i+1} + \dots + \tilde{c}_{j-1} \not\in \Z$
 for any $i < j < \epsilon(i)$.

\begin{definition}
 \label{def:4}
 Fix an admissible parameter 
 $\lambda = (\lambda_{i+1}, \dots, \lambda_{\epsilon(i)-1}) \in \C^{\epsilon(i) - i -1}$ 
 where we regard $\lambda_{\epsilon(i)} = -1$.
 We define an $\tildA_c$-module 
 $\calL_{c,\lambda}(\ii_i)$ by glueing
 local sheaves as follows:
  \begin{align*}
  \calL_{c, \lambda}(\ii_i) \vert_{X_i} &= 
  \tildA_{c, i} / \tildA_{c, i} g_i, \\
  {\calL_{c, \lambda}(\ii_i)} \vert_{X_j} &=
  \tildA_{c,j} / \tildA_{c,j} (f_j \circ g_j - \hbar \lambda_j) \\
   &=
  \tildA_{c,j} / \tildA_{c,j} (x_{\ii_j} y_{\ii_j} - \hbar \lambda_j)
  \qquad (\text{for } j = i+1, \dots, \epsilon(i)-1), \\
  {\calL_{c, \lambda}(\ii_i)} \vert_{X_{\epsilon(i)}} &=
  \tildA_{c,\epsilon(i)} / \tildA_{c,\epsilon(i)} f_{\epsilon(i)}, \\
  {\calL_{c, \lambda}(\ii_i)} \vert_{X_j} &= 0 
  \qquad (\text{for } j = 1, \dots, i-1, \epsilon(i)+1, \dots, l).
 \end{align*}
  Its glueing is given by 
\begin{equation}
\label{eq:5} 
  u_j = f_{j}^{\lambda_j - \lambda_{j+1} - \tilde{c}_j} 
 u_{j+1} \quad \text{on } X_j \cap X_{j+1}
\end{equation}
 where $u_j$ is the image of the constant function $1 \in \tildA_{c,i}$
 in $\calL_{c, \lambda}(\ii_i) \vert_{X_j}$ for $j = i$, $\dots$, $\epsilon(i)$.
\end{definition}

\begin{remark}
If $\tilde{c}_i + \tilde{c}_{i+1} + \dots + \tilde{c}_{j-1}$ for any
$j = i+1$, $\dots$, $l$, we regard $\epsilon(i) = l+1$ and the definition of
$\calL_{c,\lambda}(\ii_i)$ is given by
  \begin{align*}
  \calL_{c, \lambda}(\ii_i) \vert_{X_i} &= 
  \tildA_{c, i} / \tildA_{c, i} g_i, \\
  {\calL_{c, \lambda}(\ii_i)} \vert_{X_j} &=
  \tildA_{c,j} / \tildA_{c,j} (f_j \circ g_j - \hbar \lambda_j)
  \qquad (\text{for } j = i+1, \dots, l), \\
  {\calL_{c, \lambda}(\ii_i)} \vert_{X_j} &= 0 
  \qquad (\text{for } j = 1, \dots, i-1).
 \end{align*}
\end{remark}

Note that we have an isomorphism of $\tildA_{c,j}$-modules
\begin{equation}
 \label{eq:14}
 \calL_{c,\lambda}(\ii_i) \simeq \scM_{\lambda_j}, \qquad
 u_j \mapsto v_{\lambda_j}
\end{equation}
for $j = i+1$, $\dots$, $\epsilon(i)-1$,
under the isomorphism $\tildA_{c,j} \simeq \W_{T^* \C^1}$.

The following lemmas are proved similarly to 
\reflemma{lemma:9} and \reflemma{lemma:10}
by using \refeq{eq:14} instead of \refeq{eq:12}.

\begin{lemma}
 \label{lemma:6}
 The module $\calL_{c,\lambda}(\ii_i)$ is a well-defined 
 good $\tildA_c$-module supported on the Lagrangian subvariety
 $D_i \cup D_{i+1} \cup \dots \cup D_{\epsilon(i)}$.
 \begin{proof}
  The well-definedness is proved similarly to \reflemma{lemma:9}.

  There exists an 
  $\tildA_c(0)$-lattice 
  of $\calL_{c, \lambda}(\ii_i)$ given by
  $\calL_{c, \lambda}(\ii_i)(0) \vert_{X_j} = \tildA_c u_j$
  for $j=i$, $\dots$, $\epsilon(i)$. 
  Thus, $\calL_{c, \lambda}(\ii_i)$ is a good $\tildA_c$-module.
 \end{proof}
\end{lemma}

\begin{lemma}
 \label{lemma:7}
 For any admissible parameters
 $\lambda, \lambda' \in \C^{\epsilon(i)-i-1}$,
 we have an isomorphism of 
 $\tildA_c$-modules 
 $\calL_{c,\lambda}(\ii_i) \simeq \calL_{c,\lambda'}(\ii_i)$.
 \begin{proof}
  Note that we have $\lambda_j \not\in \Z$ because $\lambda$ is
  admissible and 
  $\tilde{c}_i + \tilde{c}_{i+1} + \dots + \tilde{c}_j-1 \not\in \Z$
  for any $i < j < \epsilon(i)$. Thus this lemma is proved similarly
  to \reflemma{lemma:10}.
 \end{proof}
\end{lemma}

By the above lemma, the $\tildA_c$-module $\calL_{c, \lambda}(\ii_i)$
is independent of the choice of the admissible parameter 
$\lambda \in \C^{\epsilon(i)-i-1}$.

\begin{definition}
 \label{def:5}
 We denote the $\tildA_c$-module $\calL_{c, \lambda}(\ii_i)$ by
 $\calL_c(\ii_i)$.
\end{definition}

In the rest of this subsection, we show that the 
$\tildA_c$-module $\calL_c(\ii_i)$ is an irreducible module.

For $i=1$, $\dots$, $l$, 
the good $\tildA_c$-module $\calL_c(\ii_i)$ is supported on
a Lagrangian subvariety 
$D_i \cup D_{i+1} \cup \dots \cup D_{\epsilon(i)}$. Thus, $\calL_c(\ii_i)$ is 
a holonomic module.
The irreducibility of the holonomic module $\calL_c(\ii_i)$ immediately 
follows from \refprop{prop:6} and \refprop{prop:7}.

\begin{proposition}
 \label{prop:8}
 The module $\calL_c(\ii_i)$ is an irreducible $\tildA_c$-module.
 \begin{proof}
  Assume there exists a nonzero submodule $\calN$ of 
  $\calL_c(\ii_i)$. By \refprop{prop:6} and \reflemma{lemma:4}, 
  we have $\Supp \calN = D_j \cup D_{j+1} \cup \dots \cup D_{k}$
  for some $i \le j \le k \le \epsilon(i)$. Assume $j \ne i$, then
  $\calL_c(\ii_i) \vert_{X_j}$ is an 
  $\tildA_{c,j} \simeq \W_{T^* \C^1}$-module and it has a
  nontrivial $\W_{T^* \C^1}$-submodule $\calN \vert_{X_j}$
  supported on $\{x=0\}$.
  On the other hand, by the definition of $\calL_c(\ii_i)$, we have
   $\calL_c(\ii_i) \vert_{X_j} \simeq \scM_{\lambda_j}$
  and $\lambda_j \not\in \Z$. By \refprop{prop:7}, 
  $\calL_c(\ii_i) \vert_{X_j}$ is an irreducible $\W_{T^* \C^1}$-module,
  and it contradict the assumption. Thus we have $j=i$. Similarly, 
  we have $k = \epsilon(i)$. Therefore $\calN = \calL_c(\ii_i)$, and thus,
  $\calL_c(\ii_i)$ is an irreducible $\tildA_c$-module.
 \end{proof}
\end{proposition}

\begin{theorem}
 \label{thm:3}
 For $i=1$, $\dots$, $l$, we have
 \[
  \gHom_{\Mod_F^{good}(\tildA_c)}(\tildA_c, \calL_c(\ii_i)) = L_c(\ii_i)
 \]
 \begin{proof}
  By \refprop{prop:7} together with the definitions of 
  $\calM^{\Delta}_c(\ii_i)$ and $\calL_c(\ii_i)$ (\refdef{def:6}, \refdef{def:4}), 
  $\calL_c(\ii_i)$ is a quotient of $\calM^{\Delta}_c(\ii_i)$. Applying the
  equivalence of \refthm{thm:1}, the $A_c$-module
  $\gHom_{\Mod_F^{good}(\tildA_c)}(\tildA_c, \calL_c(\ii_i))$ is a
  quotient of 
  $\gHom_{\Mod_F^{good}(\tildA_c)}(\tildA_c, \calM^{\Delta}_c(\ii_i)) \simeq 
  \Delta_c(\ii_i)$. Since $\calL_c(\ii_i)$ is an irreducible $\tildA_c$-module,
  the $A_c$-module $\gHom_{\Mod_F^{good}(\tildA_c)}(\tildA_c, \calL_c(\ii_i))$
  is an irreducible quotient of $\Delta_c(\ii_i)$. Therefore, it is isomorphic
  to $L_c(\ii_i)$.
 \end{proof}
\end{theorem}

Since the full subcategory of regular holonomic $\tildA_c$-modules is
closed under extensions, we have the following corollary.
\begin{corollary}
 \label{cor:2}
 For any module $M$ in $\calO(A_c)$, the corresponding $\tildA_c$-module
 $\tildA_c \otimes_{A_c} M$ is regular holonomic.
\end{corollary}

Next, we discuss the decomposition of the standard modules of $\calO(A_c)$ 
in the Grothendieck group of $\calO(A_c)$.

\begin{corollary}
 \label{cor:1}
 In the Grothendieck group of $\calO(A_c)$, we have
 \[
  [\Delta_c(\ii_i)] = \sum_{j: \tilde{c}_i + \dots + \tilde{c}_{j-1} \in \Z}
 [L_c(\ii_j)].
 \]
 \begin{proof}
  By \refprop{prop:3} (3), we have
  \[
   [\Delta_c(\ii_i)] = \sum_{j=1, \dots, l} n_j [L_c(\ii_j)]
  \]
  for some $n_j \in \Z_{\ge 0}$. 
  If $\Supp \calL_c(\ii_j) \not\subset \Supp \calM_c^{\Delta}(\ii_i)
  = D_i \cup \dots \cup D_l$, we have $n_j = 0$.
  Since $\calM_c^{\Delta}(\ii_i)$ and
  $\calL_c(\ii_j)$ are (at most) multiplicity-one on $D_k$, we have
  \[
   \sum_{j: \Supp \calL_c(\ii_j) \cap D_k} n_j = 1 \quad 
  \text{for } k = i, i+1, \dots, l.
  \]
  That is, $[\Delta_c(\ii_i)]$ is multiplicity-free in the Grothendieck
  group.
  Since $\calL_c(\ii_j)$ is a unique irreducible module whose 
  support is of the form $D_{j} \cup D_{j+1} \cup \cdots$, 
  we have $n_j = 1$ for $j = i$, $\dots$, $l$ such that
  $\tilde{c}_i + \dots + \tilde{c}_{j-1} \in \Z$ by comparing 
  the support of $\calM_c^{\Delta}(\ii_i)$ and $\calL_c(\ii_i)$.
 \end{proof}
\end{corollary}

\begin{remark}
 We can also determine the multiplicity $[\Delta_c(\ii_i) : L_c(\ii_j)]$
 in the Grothendieck group of $\calO(A_c)$
 algebraically in this case. The same result of \refcor{cor:1} 
 is immediately follows from \cite[Lemma 4.3]{K}.
\end{remark}

Finally, we discuss the subcategory of $\Mod_F^{good}(\tildA_c)$ 
corresponding to the category $\calO(A_c)$. Since a section $f$ of the
W-algebra $\tildA_c$ is invertible if and only if its symbol $\sigma_0(f)$ is
invertible in $\calO_X$, 
$\hbar^{-1} y_1 \cdots y_l$ acts locally nilpotently on an $A_c$-module $M$
if and only if $\Supp \tildA_c \otimes_{A_c} M \subset \bigcup_{i=1}^{l} D_i$. Thus, as mentioned by \cite[Remark~8.8.2]{Mc}, we have 
an equivalence of these subcategories:
\[
 \calO(A_c) \simeq \Mod_{F,\,\bigcup_{i=1}^{l} D_i}^{good}(\tildA_c)
\]
where $\Mod_{F,\,\bigcup_{i=1}^{l} D_i}^{good}(\tildA_c)$ is a
full subcategory of $\Mod_{F}^{good}(\tildA_c)$ whose modules are
supported on $\bigcup_{i=1}^{l} D_i$. As a corollary of \refcor{cor:2}, 
good $\tildA_c$-modules with $F$-action supported on $\bigcup_{i=1}^{l} D_i$
are automatically regular holonomic.

\appendix
\section{Global sections of modules}
\label{sec:glob-sect-modul}

We can calculate explicitly the global sections of 
$\calM^{\Delta}(\ii_i)$. Fix $\lambda = (\lambda_j)_{j=i+1, \dots, l} \in \C^{l-i}$
be an admissible parameter. 
First, the restriction homomorphisms
are given explicitly as follows,
\begin{align*}
 \Res_1 : \Gamma(X_j, \calM_{c, \lambda}(\ii_i)) &\longrightarrow
 \Gamma(X_j \cap X_{j+1}, \calM_{c, \lambda}(\ii_i)), \\
 f_j^m u_j &\mapsto f_j^m u_j 
 \qquad (m \in \Z_{\ge 0}), \nonumber\\
 g_j^m u_j &\mapsto C'_{-m,j} f_j^{-m} u_j
 \qquad (m \in \Z_{> 0}), \nonumber \\
 \Res_2 : \Gamma(X_{j+1}, \calM_{c, \lambda}(\ii_i)) &\longrightarrow
 \Gamma(X_j \cap X_{j+1}, \calM_{c, \lambda}(\ii_i)), \\
 g_{j+1}^m u_{j+1} &\mapsto 
 f_j^{-m + \lambda_{j+1} + \tilde{c}_j - \lambda_j} 
 u_j \qquad (m \in \Z_{\ge 0}) \nonumber \\
 f_{j+1}^m u_{j+1} &\mapsto C_{m,j+1} 
 f_j^{m + \lambda_{j+1} + \tilde{c}_j - \lambda_j} u_j,
 \qquad (m \in \Z_{\ge 0}) \nonumber \\
\end{align*}
where
\begin{align*}
 C_{m,j} &= \hbar^m (m+\lambda_j) (m+\lambda_{j}-1) \cdots (\lambda_j+1),
 \qquad (m \in \Z_{\ge 0})¡¡\\
 C'_{m,j} &= \hbar^{-m} (m+\lambda_j+1) (m+\lambda_j+2) \cdots \lambda_{j}, 
 \qquad (m \in \Z_{< 0}) \nonumber
\end{align*}
are scalar constants. For $j=i$, $\dots$, $l$ such that 
$\tilde{c}_i + \dots + \tilde{c}_{j-1} \not\in \Z$, we have
$C_{m,j}$, $C'_{m,j} \ne 0$ for all $m$.

Assume 
$\lambda = (\lambda_j)_{i+1, \dots, l} \in (\C \backslash \Z_{\ge 0})^{l-i}$.
For $j=i$, $\dots$, $l$
such that $\tilde{c}_i + \dots + \tilde{c}_{j-1} \in \Z$, we have
\begin{align}
\label{eq:19}
  C_{m,j} &\ne 0 \quad (\text{for } m < - \lambda_j, 
\text{and } j = i, \dots, l-1), \\
  C'_{m,j} &\ne 0 \quad (\text{for any } m, 
\text{and } j = i, \dots, l-1). \nonumber
\end{align}

Now, we construct the global sections of $\calM^{\Delta}_c(\ii_i)$ explicitly.
Fix $i = 1$, $\dots$, $l$ and $\lambda_{i+1}$, $\dots$, $\lambda_l$ 
such that $\lambda_j < - \tilde{c}_{i} - \tilde{c}_{i+1} - \dots - \tilde{c}_{j-1}$
for all $j=i+1$, $\dots$, $l$. 

For $j=i$, $\dots$, $l$ and $k=j$, $\dots$, $l$, set 
\[
 m_{j,k} = - \lambda_k - \tilde{c}_{k-1} - \tilde{c}_{k-2} - \dots - \tilde{c}_j.
\]
Note that we have $m_{j,k} + \lambda_{k} + \tilde{c}_{k-1} - \lambda_{k-1} =
m_{j,k-1}$.
For $j=i$, $\dots$, $l$ such that $\tilde{c}_{i} + \dots + \tilde{c}_{j-1} \in \Z$,
take $m \in \Z$ such that $0 \le m < \tilde{c}_j + \dots + \tilde{c}_{\epsilon(j)-1}$.
\footnote{we regard $\tilde{c}_l = \infty$ here.}
Then we define a section 
\begin{equation*}
 v_{j,m} = 
 \begin{cases}
  (\hbar^{-l/2} f_l)^{m_{j,l}+m} u_l
  &\text{on } X_l \\
  \left(\prod_{k=j'+1}^{l} C_{m_{j,k}+m, k} \right) 
  (\hbar^{-j'+l/2} f_{j'})^{m_{j,j'}+m} u_{j'}
  &\text{on } X_{j'} \quad (j \le j' \le l) \\
  0 &\text{on } X_{j'} \quad (j' \le j-1)
 \end{cases}
\end{equation*}
Note that 
$C_{m_{j,k}+m, k} \ne 0$ 
by \refeq{eq:19}, and $v_{j,m}$ is a well-defined global section.
Moreover, because $v_{j,m}$ is an F-equivariant section, we can identify
it with an F-equivariant homomorphism 
$\tildA_c \ni 1 \mapsto v_{j,m} \in \calM_c^{\Delta}(\ii_i)$ in
$\gHom_{\Mod_F^{good}(\tildA_c)}(\tildA_c, \calM_c^{\Delta}(\ii_i))$.

\end{document}